\documentclass[oneside,openany,a4paper]{amsart}
\DeclareFontFamily{U}{wncy}{}
    \DeclareFontShape{U}{wncy}{m}{n}{<->wncyr10}{}
    \DeclareSymbolFont{mcy}{U}{wncy}{m}{n}
    \DeclareMathSymbol{\Sh}{\mathord}{mcy}{"58} 
\usepackage[a4paper, margin=1in]{geometry}
\usepackage{amssymb}
\usepackage[centertags]{amsmath}
\usepackage{amsthm}
\usepackage{amsfonts}
\usepackage{eucal}
\usepackage{mathrsfs}
\usepackage{diagrams}

\usepackage[all,cmtip]{xy}
\usepackage{graphicx, eepic}

\usepackage{lmodern}
\usepackage[utf8]{inputenc}
\usepackage[T1]{fontenc}
\usepackage{xcolor}

\usepackage{bbm}
\usepackage{pdfpages}
\usepackage{titletoc}
\usepackage{booktabs}
\usepackage[shortlabels]{enumitem}
\usepackage{mathtools}
\usepackage{thmtools}

\usepackage{epstopdf}
\usepackage{epsfig}
\usepackage{microtype}

\usepackage{hyperref}
\hypersetup{colorlinks}
\hypersetup{
    bookmarks=true,         % show bookmarks bar?
    unicode=false,          % non-Latin characters in Acrobat’s bookmarks
    pdftoolbar=true,        % show Acrobat’s toolbar?
    pdfmenubar=true,        % show Acrobat’s menu?
    pdffitwindow=false,     % window fit to page when opened
    pdfstartview={FitH},    % fits the width of the page to the window
    pdftitle={My title},    % title
    pdfauthor={Author},     % author
    pdfsubject={Subject},   % subject of the document
    pdfcreator={Creator},   % creator of the document
    pdfproducer={Producer}, % producer of the document
    pdfkeywords={keyword1} {key2} {key3}, % list of keywords
    pdfnewwindow=true,      % links in new PDF window
    colorlinks=true,       % false: boxed links; true: colored links
    linkcolor=red,          % color of internal links (change box color with linkbordercolor)
    citecolor=magenta,        % color of links to bibliography
    filecolor=violet,      % color of file links
    urlcolor=blue      % color of external links
}
\usepackage{multirow, url}
\usepackage{textcomp}
\DeclareMathAlphabet{\cmcal}{OMS}{cmsy}{m}{n}
\newtheoremstyle{thm}% name of the style to be used
  {3pt}% measure of space to leave above the theorem. E.g.: 3pt
  {3pt}% measure of space to leave below the theorem. E.g.: 3pt
  {\em}% name of font to use in the body of the theorem
  {0pt}% measure of space to indent
  {\bfseries}% name of head font
  {}% punctuation between head and body
  {5pt}% space after theorem head
  {}% Manually specify head
\newtheoremstyle{rem}% name of the style to be used
  {3pt}% measure of space to leave above the theorem. E.g.: 3pt
  {3pt}% measure of space to leave below the theorem. E.g.: 3pt
  {}% name of font to use in the body of the theorem
  {0pt}% measure of space to indent
  {\bfseries}% name of head font
  {.}% punctuation between head and body
  {5pt}% space after theorem head
  {}% Manually specify head

\newtheorem{thm}{Theorem}[section]

\newtheorem{lem}[thm]{Lemma}

\newtheorem{conj}[thm]{Conjecture}

\theoremstyle{definition}

\theoremstyle{rem}

\numberwithin{equation}{section} \numberwithin{table}{section}
	
\newtheorem*{thm*}{Theorem}
\newtheorem*{rem*}{Remark}
\newtheorem*{rems*}{Remarks}
\newtheorem*{exam*}{Example}
\newtheorem*{exams*}{Examples}

\makeatletter
\newcommand{\neutralize}[1]{\expandafter\let\csname c@#1\endcsname\count@}
\makeatother

\newcommand{\C}{{\mathbb{C}}}

\newcommand{\F}{{\mathbb{F}}}

\newcommand{\Q}{{\mathbb{Q}}}
\newcommand{\R}{{\mathbb{R}}}

\newcommand{\bV}{{\mathbb{V}}}

\newcommand{\Z}{{\mathbb{Z}}}

\newcommand{\fb}{{\mathfrak{b}}}

\newcommand{\cA}{{\cmcal{A}}}

\newcommand{\cO}{{\cmcal{O}}}
\newcommand{\cP}{{\cmcal{P}}}

\newcommand{\cU}{{\cmcal{U}}}

\def\n2Z{\frac{1}{n^2}\Z/\Z}

\def\z{\zeta}

\def\bq{\begin{quote}}
\def\eq{\end{quote}}

\newcommand{\Gm}{{\mathbb{G}}_{m}}

\newcommand{\zmod}[1]{{\Z/{#1}\Z}}

\newcommand{\arinj}{\ar@{^(->}}
\newcommand{\arsurj}{\ar@{->>}}
\newcommand{\arsub}{\ar@{}[r]|-*[@]{\subset}}
\newcommand{\arsup}{\ar@{}[r]|-*[@]{\supset}}
\newcommand{\arcap}{\ar@{}[d]|-*[@]{\subset}}
\newcommand{\arcup}{\ar@{}[u]|-*[@]{\subset}}
\newcommand{\arin}{\ar@{}[u]|-*[@]{\in}}

\renewcommand{\mod}[1]{{\,\operatorname{mod}\hspace{0.5mm} {#1}}}

\usepackage{cleveref}

\newcommand{\Ker}{\operatorname{Ker}}

\newcommand{\Hom}{{\operatorname{Hom}}}
\newcommand{\Gal}{{\operatorname{Gal}}}

\newcommand{\Aut}{{\operatorname{Aut}}}

\newcommand{\Spec}{{\operatorname{Spec}}}

\newcommand{\inv}{{\operatorname{inv}}}

\newcommand{\Tr}{{\operatorname{Tr}}}

\renewcommand{~}{\hspace{0.4mm}}

\newcommand{\ov}{\overline}

\newcommand{\bd}{\boldsymbol}

\mathchardef\hyp="2D

\def\rank{\mbox{rank}}

\newcommand{\xyv}[1]{\xymatrixrowsep{#1 pc}}

\def\Kbar{\bar{K}}

\def\Qbar{\overline{\mathbb{Q}}}
\def\hZ{\widehat{\mathbb{Z}}}
\def\loc{\mbox{loc}}

\def\cU{\mathcal{U}}

\def\<{\langle }
\def\>{\rangle}
\def\hQp{\widehat{\Q_p^*}}
\def\hF{\widehat{F^*}}

\def\bF{\bar{F}}
\def\nZ{\frac{1}{n}\mathbb{Z}/\mathbb{Z}}

\def\ms{\medskip}

\def\invlim{\varprojlim}

\def\A{\mathbb{A}}

\def\bd{\begin{diagram}}
\def\ed{\end{diagram}}
\def\bV{\overline{V}}
\def\tc{\tilde{c}}
\renewcommand{\O}{{\cmcal O}}

\title{Arithmetic Gauge Theory: A Brief Introduction}
\author{Minhyong Kim}
\address{Mathematical Institute, University of Oxford, Woodstock Road, Oxford, OX2 6GG, and The Korea Institute for Advanced Study 85 Hoegiro, Dongdaemungu,
Seoul 02455,
Republic of Korea}
\thanks{Supported by  grant 	EP/M024830/1 from the EPSRC}
\begin{document}

\maketitle
\begin{abstract}
Much of  arithmetic geometry is concerned with the study of principal bundles. They occur prominently in the arithmetic of elliptic curves and, more recently, in the study of the Diophantine geometry of curves of higher genus. In particular, the geometry of moduli spaces of principal bundles appears to be closely related to an effective version of Faltings's theorem on finiteness of rational points on curves of genus at least 2. The study of arithmetic principal bundles includes the study of {\em Galois representations}, the structures linking motives to automorphic forms according to the Langlands programme. In this article, we give a brief introduction to the arithmetic geometry of principal bundles with emphasis on some elementary analogies between arithmetic moduli spaces and the constructions of quantum field theory. For the most part, it can be read as an attempt to explain standard constructions of arithmetic geometry using  the language of physics, albeit employed in an amateurish and ad hoc manner. \end{abstract}

\section{Fermat's principle}Fermat's principle says that the trajectory taken by a beam of light is a solution to an optimisation
problem. That is, among all the possible paths that light could take, it selects the one requiring the least time to traverse. This was the first example
of a very general methodology known nowadays as the principle of least action. To figure out the trajectory or spacetime configuration favoured by nature, you
should analyse the physical properties of the system to associate to each possible configuration a number,
called the action of the configuration. Then the true trajectory is one where the action is extremised. The action determines a constraint equation, the so-called Euler-Lagrange equation of the system,
whose solutions give you possible trajectories. The action principle in suitably general form is the basis of classical field theory,
particle physics, string theory, and gravity.  For Fermat to
have discovered this idea so long ago in relation to the motion of light was a monumental achievement,
central to the scientific revolution that rose out of  the intellectual fervour of 17th century Europe.

However, Fermat is probably better known these days as the first modern number-theorist. Among the
intellectual giants of the period,  Fermat was
almost unique in his preoccupation with prime numbers and Diophantine equations, polynomial equations
to which one seeks integral or rational solutions. Located among his many forays into this subject one
finds his famous `Last Theorem', which elicited from the best mathematical minds of subsequent generations
several hundred years of theoretical development before it was finally given a proof by Andrew Wiles in
1995 \cite{wiles}. The action principle and Fermat's last theorem are lasting tributes to one of the singularly original
minds active at the dawn of modern science. Could there be a relation between the two? In fact, the problem of finding the trajectory of light and that of finding rational solutions to Diophantine equations are two facets of the same problem, one occurring in geometric gauge theory, and the other, in {\em arithmetic gauge theory}. The fact that the photon is described by a $U(1)$ gauge field is well-known. The purpose of this article is to give the motivated physicist with  background in geometry and topology some sense of the second type of theory and its relevance to the theory of Diophantine equations. 

In the context of abelian problems, say the arithmetic of elliptic curves, much of the material is classical. However, for non-abelian gauge groups,  the perspective of gauge theory is very useful and has concrete consequences. It may be that number-theorists  can also benefit from the intuition provided by this somewhat fanciful view, even though much of it will appear as pretentious reformulation of well-known notions.  In any case, it is hoped that the impressionistic treatment of this paper will not be overly irritating, since it is mostly supplemented by pointers to published literature.\section{Diophantine geometry and gauge theory}

We will be employing the language of Diophantine geometry, whereby a system of equations is encoded in an algebraic variety 
$$V$$
defined over $\Q$. We will always assume $V$ is connected. The rational solutions (or points, to use the language of geometry) will be denoted by $V(\Q)$, while $p$-adic and adelic points\footnote{ The main advantage of the field of $p$-adic numbers   over the reals is its substantial but manageable absolute Galois group. The adeles can be thought of as essentially the product ring of $\R$ and $\Q_p$ for all $p$, with some small restriction. See \cite{neukirch} for a review.} will be denoted by $V(\Q_p)$ and $V(\A_{\Q})$, respectively.  Sometimes integral models will be implicitly assumed, in which case, we will write $V(\Z)$ for the integral points. Similarly, $V(B)$ will denote the points of $V$ in a general ring $B$. For example, one common ring occuring in Diophantine geometry is $\Z_S$, the $S$-integers for some finite set $S$ of primes, consisting of rational numbers that only have primes from $S$ in the denominator, and therefore, intermediate between $\Z$ and $\Q$. When we think of  rings as functions,  a ring like $\Z_S$ should be thought of as those having singularities lying in a fixed set.

\ms

Even though we will not use much of it, we remark without explanation the formulation  in the language of schemes, whereby a $B$-point can be viewed as a section $s$ of a fibration over $\Spec(B)$:

$$
\xy
(0,22)*{V_B};
(0,-2)*{\Spec(B)};
(0,20)*{}="A";
(0,0)*{}="B";
(3,0)*{}="C";
(3,20)*{}="D";
{\ar@{->}"A";"B"};
{\ar^{s}@/_1pc/"C";"D"};
\endxy
$$
\ms

The theory of \cite{kim1, kim2,kim3, kim4} associates to $p$-adic or adelic points of $V$
 {\em arithmetic gauge fields}. We will be focussing mostly on the $p$-adic theory for the sake of expositional simplicity. The statement, which we will review in sections 4 and 5, is that there is a natural map
 $$A_p: V(\Q_p)\rTo \mbox{$p$-adic arithmetic gauge fields}$$
 The type of gauge field is determined by the arithmetic geometry of $V$. Among $p$-adic or adelic gauge fields, the problem is to find the locus of {\em rational gauge fields}. The condition for a gauge field to be rational or integral\footnote{For the most part, our varieties will be projective, allowing us to identify integral and rational points.} can be phrased entirely in terms of global symmetry, and is shown in many cases to impose essentially computable constraints on the $p$-adic gauge fields. These constraints should be viewed as one version of `arithmetic Euler-Lagrange (E-L) equations'. In number theory, they are closely related to {\em reciprocity laws} as will be explained in section 8 and section 10. The key point is that when the solution $x\in V(\Q_p)$ lies in the subset $V(\Q)$, then the corresponding gauge field $A_p(x)$ is rational.  That is, we have a commutative diagram
 $$\bd X(\Q)&\rInto & X(\Q_p)\\
 \dTo^A & & \dTo^{A_p} \\
 \mbox{rational gauge fields}&\rInto &\mbox{$p$-adic gauge fields}.\ed$$
 \ms
 
 The E-L equation for $A_p(x)$ can be translated,  using $p$-adic Hodge theory,  back to an analytic equation satisfied by the  point $x$. When $V$ is a curve and the equation thus obtained is non-trivial, this implies finiteness theorems for rational points. That is, it is often possible to prove that $$A_p^{-1}(\mbox{rational gauge fields})$$ is a finite set\footnote{In fact, this will always be true subject to standard conjectures on mixed motives \cite{kim6}.}.
One can  give thereby new proofs of the finiteness of rational solutions to a range of Diophantine equations, including the generalised Fermat equations \cite{CK}
$$ax^n+by^n=c$$
for $n\geq 4$.
This finiteness was first proved by Gerd Faltings in 1983 as part of his proof of the Mordell conjecture (cf. Section 3) using ideas and constructions of arithmetic geometry,
 However, the  proof in  \cite{CK} has a number of theoretical
advantages as well as practical ones. On the one hand, the gauge-theoretical perspective has the potential to be applicable  to a very broad class of phenomena encompassing
many of the central problems of current day number theory  \cite{kimcs}. On the other, unlike Faltings's proof, which is
widely regarded as ineffective, the gauge-theory proof conjecturally leads to a computational method for actually finding rational
solutions \cite{kim5}, a theme that is currently under active investigation \cite{BD1, BD2, DW1, DW2}.

It should be remarked that the map $A$ that associates gauge fields to points has been well-known since the 1950s when the variety $V$ is an elliptic curve, an abelian variety, or generally, a commutative algebraic group. More general equations, for example, curves of genus $\geq 2$, require non-abelian gauge groups, and  it is in this context that the analogy with physics assumes greater important. Nonetheless, the arithmetic E-L equations obtained thus far have not been entirely canonical. The situation is roughly that of having an Euler-Lagrange equation without an action. On the other hand, if we consider gauge theory with constant gauge groups (to be discussed below), there is a very natural analogue of the Chern-Simons action on  3-manifolds, for which it appears a theory can be developed in a manner entirely parallel to usual topology. In particular, some rudiment of path integral quantisation becomes available, and give interpretations of $n$-th power residue symbols as arithmetic linking numbers \cite{kimcs, CKKPY, CKKPPY}. The arithmetic Chern-Simons action of those papers were originally motivated by the problem of defining an action for gauge fields arising in Diophantine geometry.

\section{ Principal bundles and number theory: Weil's constructions}
In the language of geometry, gauge fields are {\em  principal bundles with connection}, and this is the form in which we will be discussing arithmetic analogues. Perhaps it is useful to recall that over the last 40 year or so, the  idea that a space $X$ can be fruitfully studied in terms of the field theories it can support has been extraordinarily powerful in geometry and topology. The space of interest can start out both as a target space of fields or as a  source. Both cases are able to give rise to suitable moduli spaces of principal bundles (with connections) $$M(X, G),$$ which then 
can be viewed as invariants of $X$.  Here, $G$ might be a compact Lie group or an algebraic group, while the moduli space might consist of flat connections, or other spaces of solutions  to differential equations, for example, the (self-dual) Yang-Mills equation. Of course this idea is at least as old as Hodge theory for abelian $G$, while the non-abelian case has seen an increasing  array of deep interactions with physics since the work of Atiyah, Bott, Drinfeld, Hitchin, Manin, Donaldson, Simpson, Witten, and many others \cite{atiyah, AHDM, AB, donaldson, simpson, witten}.
\ms

However,
my impression is that it is not widely known among mathematicians that the study of principal  bundles  was from its inception closely tied to number theory. Probably, the first moduli space of principal bundles   appeared  as the  Jacobian  of a Riemann surface,  the complex torus target of the Abel-Jacobi map
$$x\mapsto (\int_b^x \omega_1,\int_b^x \omega_2, \ldots, \int_b^x \omega_g) \mod H_1(X, \Z),$$
in what can now be interpreted as the Hodge realisation. 
 In the early 20th century \cite{weil1},  Andre Weil gave the first algebraic construction of the Jacobian $J_C$ of a smooth projective algebraic curve $C$ of genus at least two defined over an algebraic number field $F$. His main motivation was the {\em Mordell conjecture}, which said that the set $C(F)$ of $F$-rational points should be finite. The algebraic nature of the construction allowed $J_C$ to be viewed also as a variety in its own right over $F$ that admitted an embedding
 $$C(F)\hookrightarrow J_C(F)$$
 $$x\mapsto \cO(x)\otimes \cO(-b).$$
 Weil proved that $J_C(F)$ was a finitely-generated abelian group, but was unable to use this striking fact to prove the finiteness of $C(F)$. It appears to have taken him another decade or so \cite{weil2} to realise that the abelian nature of $J_C$ kept it from being too informative about $C(F)$, and from there he went on to define
 $$  M(C, GL_n)(F) = [\prod_{x\in C} GL_n(\cO_{C,x})]\backslash GL_n(\A_{F(C)})/GL_n(F(C))$$
the set of isomorphism classes of rank $n$ vector bundles on $C$. Weil considered this as a non-abelian extension of the Jacobian, which might be applied to the non-abelian arithmetic of $C$. Even though the Mordell conjecture remained unproven for another 45 years, Weil's construction went on to inspire many ideas in geometric invariant theory and non-abelian Hodge theory, much of it in interaction with Yang-Mills theory \cite{AB, NS, donaldson, simpson}.
\ms

In order to bring about further applications to number theory, it  turned  out to be critical to consider moduli of principal bundles over $F$ itself, or over various rings of integers in $F$, not just over other objects of algebro-geometric nature sitting over $F$. These are the arithmetic gauge fields mentioned above.

\section{ Arithmetic gauge gields}
For the most part, in this paper, we will present the theory in a pragmatic manner, requiring as little theory as possible. What underlies the discussion is the topology of the spectra of number fields, local fields, and rings of integers, but it is possible to formulate most statements in the language of fields and groups. Roughly speaking, when we refer to an object over (or on) a ring $\mathcal{O}$, we will actually have in mind the geometry $\Spec(\mathcal{O})$, the spectrum\footnote{But the reader will not be required to know the language of spectra or schemes until section 9.} of $\mathcal{O}$.
\ms

 Given a field $K$ of characteristic zero, denote by $$G_K=\Gal(\Kbar/K)$$ the Galois group of an algebraic closure $\Kbar$ of $K$. Thus, these are the field automorphisms of $\Kbar$ that act as the identity on $K$. For any finite extension $L$ of $K$ contained in $\Kbar$, the algebraic closure $\bar{L}$ is the same as $\Kbar$,  and
 $$G_L=\{g\in G_K \ | \ g|L=I\}.$$ When $L/K$ is Galois, $G_L$ is the kernel of the projection $$G_K\rTo \Gal(L/K).$$ In fact, we can write the Galois group   as an inverse limit \footnote{An element of such an inverse limit is a compatible collection $(g_L)_L$,
 where the compatibility means that if $L\supset L'\supset K$, then $g_L|L'=g_{L'}$.}
 $$G_K=\invlim_L \Gal(L/K),$$
 as $L$ runs over the finite Galois extensions of $K$ contained in $\Kbar$.  This equips $G_K$ with the topology of a compact, Hausdorff, totally disconnected space, with a basis of open sets given by the cosets of the $G_L$.  In particular, it is homeomorphic to   a Cantor set.
Such large inverse limits afford some initial psychological difficulty, but form an essential part of arithmetic topology. 
 
 By {\em a gauge group over $K$} we mean a topological group $U$ with a continuous action of $G_K$.  A {\em $U$-gauge field,} or {\em principal $U$-bundle} over $K$ is a topological space $P$ with a simply-transitive continuous right $U$-action and a continuous left $G_K$ action that are compatible. This means
$$g (p u)=g(p)g(u)$$
for all $g \in G_K$, $u\in U$ and $p\in P$. 
We remark that a principal $G$-bundle as defined  corresponds naively only to flat connections in geometry. We will comment on this analogy in more detail below. In algebraic geometry, the expression {\em $U$-torsor} is commonly used in place of the differential geometric terminology. We will use both. For the purposes of this paper,  an {\em arithmetic gauge group (or field)}  will mean a gauge group (or field) over an algebraic number field or a completion of an algebraic number field.

There is an obvious notion of isomorphism of $U$-torsors, and a well-known classification of $U$-torsors over $K$: 
Given $P$, choose $p\in P$. Then for any $g\in G_K$, $g(p)=pc(g)$ for a unique $c(g)\in G_K$. It is easy to check that  $g \mapsto c(g)$ defines a continuous function $$c:G_K\rTo U$$ such that
$$c(gg')=c(g)gc(g').$$
The set of such functions is denoted $Z^1(G_K, U)$, and called the set of continuous 1-cocycles of $G_K$ with values in $U$. There is a right action of $U$ on $Z^1(G_K, U)$ by
$$(uc)(g)=g(u^{-1})c(g)u,$$
and  we define
$$H^1(G_K, U):=Z^1(G_K,U)/U.$$

\begin{lem}
The procedure described above defines a bijection
$$\mbox{Isomorphism classes of $U$-torsors}\simeq H^1(G_K, U).$$
\end{lem}

We will denote $H^1(G_K,U)$ also by $H^1(K, U)$, to emphasise its dependence on the topology of $\Spec(K)$.
A rather classical case is when $U=R(\Kbar)$, the $\Kbar$-points of an algebraic group  $R$ over $K$, which we consider with the discrete topology. 
We will often write $R$ for $R(\Kbar)$, when there is no danger of confusion.  A trivial but important example is  $R=\Gm$, the multiplicative group, so that  $\Gm(\Kbar)=\Kbar^{\times}$. In this case, Hilbert's theorem 90 \cite{NSW} says
$$H^1(K, \Gm)=0,$$
or that every principal $\Gm$-bundle is trivial.
Another important class  is  that of abelian varieties for example,  elliptic curves. In that case $H^1(K, R)$ is usually called the Weil-Chatelet group of $R$ \cite{silverman}. 

Some useful operations on torsors include
\ms

(1) Pushout: If $f: U\rTo U'$ is a continuous homomorphism of groups over $K$, then there is a pushout functor $f_*$ that takes $U$-torsors to principal $U'$-torsors. The formula is
$$f_*(P)=[P\times U']/U,$$
where the right action of $U$ on the product is $(p,u')u=(pu, f(u^{-1})u').$ The resulting quotient still has the $U'$-action: $[(p,u')]v=[(p, u'v)].$
\ms

(2) Product: When $P$ is an $U$-torsor and $P'$ is an $U'$-torsor, $P\times P'$ is a $U\times U'$-torsor.

\ms

Note that if $U$ is abelian, the group law $$m:U\times U\rTo U$$ is a  homomorphism. Using this, 
$$(P,P')\mapsto m_*(P\times P')$$
defines a bifunctor on principal $U$ bundles and an abelian group law on $H^1(K, U)$. However, if $U$ is non-abelian, there is no group structure on the $H^1$ and  matters becomes more subtle and interesting.

When $U$ is an abelian group, one can define cohomology groups in every degree
$$H^i(K, U):=\Ker[d: C^i(G_K, U)\rTo C^{i+1}(G_K, U)]/ Im[ d: C^{i-1}(G_K, U)\rTo C^{i}(G_K, U)].$$
Here, $C^i(G, U)$ is the set of continuous maps from $G^i$ to $U$, while the differential $d$ is defined in a natural combinatorial manner \cite{NSW}. One checks that
$H^0(K,U)=U^{G_K}$, the set of invariants of the action, and that the cohomology groups fit into a long exact sequence as usual. That is, if
$$1\rTo U''\rTo U\rTo U'\rTo 1$$
is exact, then we get
$$0\rTo (U'')^{G_K}\rTo U^{G_K}\rTo (U')^{G_K}\rTo H^1(K, U'')\rTo H^1(K, U)\rTo H^1(K,U') $$
$$\rTo H^2(K, U'')\rTo H^2(K,U)\rTo H^2(K,U')\rTo \cdots$$
The sequence up to the $H^1$ terms remains exact even when the groups are non-abelian, except the meaning needs to be interpreted a bit carefully.

An important case is  $U=R(\Kbar)$ for $R$ a connected abelian algebraic group. In this case, multiplication by $n$ induces an exact sequence
$$0\rTo R[n]\rTo R\rTo^n R\rTo 0,$$
where $A[n]$ generally denotes the $n-$torsion subgroup of an abelian group $A$. Hence, we get the long exact sequence
$$0\rTo (R[n])^{G_K}\rTo R^{G_K}\rTo R^{G_K}\rTo H^1(K, R[n])\rTo H^1(K, R)\rTo H^1(K,R)\rTo $$
Note here that $R^{G_K}=R(K)$, the $K$-rational points of $R$. Thus, we get an injection
$$R(K)/nR(K)\rInto H^1(K, R[n]),$$
indicating how principal bundles for $R[n]$ can encode information about the group of rational points. When $R$ is an elliptic curve, this is the basis of the descent algorithm for computing the Mordell-Weil group, about which we will say more later.

Some genuinely topological groups $U$ arise from taking inverse limits. For example, we have the group $\mu_n\subset \Gm$ of $n$-th roots of unity.
They are related by the system of power maps
$$\mu_{ab}\rTo^{(\cdot)^a}\mu_b,$$
so that we can take an inverse limit
$$\hZ(1):=\invlim_n \mu_n.$$
This is a topological group isomorphic to $\hZ$, the profinite completion\footnote{ Given any group $A$, the profinite completion $\hat{A}$ of $A$ is by definition
$$\hat{A}=\invlim_N A/N,$$
where $N$ are normal subgroups of finite index.} of $\Z$, but with a non-trivial action of $G_K$. 
It is common to focus on a set of prime powers for a fixed prime $p$, and define
$$\Z_p(1)=\invlim \mu_{p^n}.$$
As a topological group, $\Z_p(1)\simeq \Z_p$, the group of $p$-adic integers\footnote{Recall that the $p$-adic integers are sometime represented as power series $\sum_{i=0}^{\infty} a_ip^i$ with $0\leq a_i\leq p-1$. Another representation is $\Z_p=\invlim \Z/p^n$.}. It is a simple example of a compact $p$-adic Lie group.
Principal bundles for this are then classified by $H^1(K, \Z_p(1))$.

In fact, we have an isomorphism
$$H^1(K, \hZ(1))\simeq \invlim H^1(K, \mu_{n}),$$
so we can consider a $ \hZ(1))$-torsor as being a compatible collection of $\mu_{n}$-torsors as we run over $n$. The exact sequence
$$1\rTo \mu_{n}\rTo \Gm\rTo^{n}\Gm\rTo 1$$
gives rise to the long exact sequence
$$1\rTo \mu_{n}(K)\rTo K^{\times}\rTo^nK^{\times} \rTo H^1(K, \mu_{n})\rTo 0,$$
so that we get an isomorphism
$$K^{\times}/(K^{\times})^{n}\simeq H^1(K,\mu_n).$$
Concretely, the torsor associated to an element $a\in K$ is simply the set $a^{1/n}$ of $n-$th roots of $a$ in $\Kbar$. This clearly admits an action
of $\mu_{n}$. That is, the group $\mu_n$ can be thought of as `internal symmetries' of the set $a^{1/n}$. This torsor  only depends on the class of $a$-modulo $n$-th powers, and is trivial if and only if $a$ has an $n-$th root in $K$. The point is that the choice of any $n$-th root in $\Kbar$ will determine a bijection to $\mu_n$, but this will be equivariant for the $G_K$-action exactly if you choose an $n$-th root in $K$ itself, which may or may not be possible. In discussing torsors over fields, it will be important in this way to keep track of both the $U$-action, the internal symmetries, and the $G_K$-action, which can be thought of as the analogue of external (spacetime) symmetries in physics.

\section{Homotopy  and gauge fields}
We will generalise the discussion of internal and external symmetries of the previous section. Let $V$ be a variety defined over $K$ and $b\in V(K)$ a $K$-rational point\footnote{ In most of the work thus far, a basepoint $b$ was used. It is possible to develop the theory without such a choice. But then, instead of a moduli space of torsors, we will be dealing with a {\em gerbe}.}. From this data one gets a gauge group as well as torsors on $K$ associated to rational points of $V$.
The gauge group will be $$U=\pi_1(\bV,b),$$ one of the many different versions of the fundamental group\footnote{We will not give the precise definitions in terms of fibre functors. A good general introduction is the book of Szamuely \cite{szamuely}, while the algebraic group realisation we will use below is given a careful discussion in \cite{deligne}.} of $\bV$,  which is $V$ regarded\footnote{ The general principle is that varieties over algebraic closed fields belong to the realm of usual geometry, while there is always an arithmetic component to geometry over non-closed field. But even in dealing with such subtleties, one constantly uses geometry over the algebraic closure.} as a variety over $\Kbar$. We will not take care to distinguish notationally between different types of fundamental groups, since the context will make it clear which one is being referred to. (Conceptually, it is also useful to regard them all as essentially the same.)
Whenever $K$ is embedded into $\C$, it will be  a completion of the topological fundamental group of $V(\C)$, either in a profinite or an algebraic sense. However, the key point is that it admits an action of $G_K$, and has the structure of a gauge group over $K$. The $G_K$-action is usually highly non-trivial, and this is a main difference from geometric gauge theory, where the gauge group tends to be constant over spacetime. Now, given any other point $x\in V(K)$, we associate to it the homotopy classes of path $$P(x):=\pi_1(\bV;b,x)$$ from $b$ to $x$, which then has both a compatible action of $G_K$ and of $\pi_1(\bV,b)$. That is, 

\begin{quote}{\em 
the loops based at $b$ are acting as internal symmetries of sets of paths emanating from $b$, while $G_K$ acts compatibly as external symmetries\footnote{This is a very elementary idea, but worth emphasising in my view.}. }
\end{quote}
In order to provide some intuition for the $G_K$-action, we give a rather concrete description in the case where $\pi_1(\bV,b)$ is the profinite \'etale fundamental group. It is worth stressing again that this is just the profinite completion of the topological fundamental group of $V(\C)$, the complex manifold associated to $V$ via some complex embedding of $K$.  However, the remarkable, albeit elementary, fact is the existence of the `hidden' Galois symmetry.
To describe it,  one approach is to construct the fundamental group and path spaces using covering spaces. 

Recall that for a manifold $M$, if $$f: \tilde{M}\rTo M$$ is the universal covering space, then the choice of a basepoint $m\in M$ and  a lift $\tilde{m}\in \tilde{M}_m:=f^{-1}(m)$  determines a canonical bijection
$$\pi_1(M,m)\simeq \tilde{M}_m$$
that takes $e$ to $\tilde{m}$. This bijection is induced by the homotopy lifting of paths. Similarly, 
$$\pi_1(M;m, m')\simeq \tilde{M}_{m'}.$$
If we replace $M$ by the variety $\bV$, there is still a notion of an algebraic universal covering $$\tilde{\bV}\rTo \bV,$$ except it is actually an inverse system $$(\bV_i\rTo \bV)_{i\in I}$$ of finite algebraic covers, each of which is unramified, that is, surjective on tangent spaces. The universality means that any finite connected unramified cover  is dominated by one of the $\bV_i$.  For an easy example, consider the compatible system of $n$ power maps $$(\bar{\mathbb{G}}_m\rTo^{(\cdot)^n}\bar{\mathbb{G}}_m)_n.$$These together form the algebraic universal cover $\tilde{\bar{\mathbb{G}}}_m\rTo\bar{\mathbb{G}}_m$. 

Now if we choose a basepoint $b\in V(K)$ and a lift\footnote{ By this, we mean a compatible system $b_i$ of basepoint lifts to the finite covers $\bV_i\rTo \bV$.
Compatibility here means that whenever you have  a map $\bV_i\rTo \bV_j$ of covers, $b_i$ is taken to $b_j$. } $\tilde{b}\in \tilde{\bV}$, then there is a unique $K$-model $$\tilde{V}\rTo V$$ of $\tilde{\bV}$, that is, a system $$(V_i\rTo V)_i$$ defined over $K$ that gives rise to the universal covering over $\Kbar$, characterised by the property that $\tilde{b}$ consists of $K$-rational points of the system\footnote{One way to see this is that the {\em pointed} covering $(\tilde{\bV}, \tilde{b})\rTo ( \bV,b)$ is really universal, in that  any other pointed covering is dominated by a {\em unique } map from $\tilde{\bV}$.  By applying this to  Galois conjugates of $(\tilde{\bV}, \tilde{b})$, we get descent data that give a $K$-model for the pointed system. This kind of reasoning is usually called `Weil descent'. For details, see \cite{milne}.}.

Even though we have not given a formal definition of the profinite \'etale fundamental group, a useful fact is that there are canonical bijections
$$\pi_1(\bV,b)\simeq \tilde{V}_b$$
and
$$\pi_1(\bV;b,x)\simeq \tilde{V}_x.$$
That is, the fundamental group and the homotopy class of paths can be identified with the fibers of the universal covering space.  This way of presenting them makes it somewhat hard to see the torsor structure. On the other hand, it does make it apparent how $G_K$ is acting. The problem of describing this action can be thought of as that of giving some manageable construction of $\tilde{V}$. This is in general a quite hard problem and typically, one studies some quotient of the fundamental group corresponding to special families of covers, such as abelian covers or solvable covers. An alternative is to study  linearisations of the fundamental group, which we will discuss below.

Anyway, we end up with a map
$$V(K)\rTo H^1(G_K, \pi_1(\bV, b));$$
$$x\mapsto \pi_1(\bV;b,x);$$
encoding points of $V$ into torsors. Typically $H^1(G_K, \pi_1(\bV, b))$ will be much bigger than $V(K)$. That is, there will be many torsors\footnote{However, there are important cases where this is conjectured to be a bijection. This is the subject of Grothendieck's {\em section conjecture} \cite{grothendieck}.} that are not of the form $P(x)$ for some point of $x$. But the important thing for us is that the space of torsors often carries a natural geometry, remarkably similar to the geometry of classical solutions to a geometric gauge theory. This added  geometric structure turns out to be very useful in grappling with the sparse structure of $V(K)$.
\section{The local to global problem, reciprocity laws, and Euler-Lagrange equations}
For the remainder of this paper, we will assume that $U$ is either a $p$-adic Lie group\footnote{We will not define this notion here, but rely on examples like $\Z_p$, $GL_n(\Z_p)$,   $p$-adic points of more general reductive algebraic groups, finite groups, and  group extensions formed out of such groups. For a systematic treatment, see \cite{schneider}.} for a fixed prime $p$ or a discrete group. (Depending on convention, the latter can be included in the former.)
So as to avoid discussing detailed algebraic number theory \cite{CF}, we will  focus mostly on $K=\Q$ or $K=\Q_v$, where $v$ could be a prime $p$ or the symbol $\infty$. We will refer to any such $v$ as a {\em place} of $\Q$, as it corresponds to an equivalence class of absolute values. The field $\Q_v$ is obtained by completing with respect to an absolute value corresponding to $v$. Thus, we have the field $\Q_p$  of $p-$adic numbers, while $\Q_{\infty}$ denotes the field of real numbers $\R$.

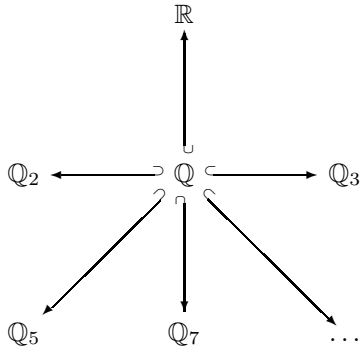
$$\bd & & \R & & \\
& & \uInto & &  \\
\Q_2  & \lInto & \Q  & \rInto &\Q_3 \\
  &\ldInto & \dInto&\rdInto & \\
 \Q_5 & & \Q_7& & \ldots 
  \ed$$

We will denote by $\Qbar$ the field of algebraic numbers and $$\pi:=G_{\Q}=\Gal(\Qbar/\Q).$$We denote by $\Qbar_v$  an algebraic closure of $\Q_v$
and $$\pi_v:=G_{\Q_v}=\Gal(\Qbar_v/\Q_v).$$ For each v, we choose an embedding $\Qbar\rInto \Qbar_v$.  Restricting the action of $\pi_v$ to $\Qbar$ then determines an embedding \footnote{ The fact that this is an embedding is not entirely obvious. It has to do with the denseness of algebraic numbers inside $\Qbar_v$. The reader should be aware that $\pi_v$ is a very thin subgroup of $G$ . It is topologically finitely generated and has an explicit description \cite{NSW}. The structure of $G$, on the other hand, is still very mysterious.}
$$\pi_v\rInto \pi$$ for each $v$. 

We will need one more mildly technical fact about the structure of $\pi_p$ for primes $p$ \cite{neukirch}. The field $\Q_p$ has an integral subring $\Z_p$, the $p$-adic integers. The integral closure\footnote{ This refers to the elements of the field extension that satisfy a nontrivial monic polynomial equation with coefficients in $\Z_p$. This notion is most natural when we consider the integral closure of $\Z$ in a field extension $F$ of $\Q$ of dimension $d$. In this setting, the integral closure is the maximal subring of $F$ isomorphic to $\Z^d$ as a group.}  of $\Z_p$ in $\Qbar_p$ is a subring $$\O_{\Qbar_p}\subset \Qbar_p$$ that  is stabilised by the $\pi_p$-action. The ring $\O_{\Qbar_p}$ has a unique maximal ideal $m_p$, and $$\O_{\Qbar_p}/m_p\simeq \bar{\F}_p,$$ an algebraic closure of $\F_p$. Thus, acting on this quotient ring determines a homomorphism
$$\pi_p\rTo \pi_p^u=\Aut(\O_{\Qbar_p}/m_p)\simeq \Gal(\bar{\F}_p/\F_p)$$
that turns out to be surjective\footnote{The superscript is supposed to stand for `unramified', corresponding to the fact that $\pi_p^ur$ is the Galois group of the  maximal extension in $\Qbar_p$ that is unramified over $\Q_p$.}. The last Galois group is generated by one element $Fr'_p$ whose effect is   $x\mapsto x^{1/p}$. (This is the group-theoretic inverse of the usual generator, the $p$-th power map\footnote{The reason for using the inverse rather than the natural $p$-power map has to do with the geometric Frobenius map acting on \'etale cohomology. This is a rather confusing convention, about which I would suggest the reader refrain from asking further at the moment.}.) Any lift $Fr_p$ of $Fr'_p$ to $\pi_p\subset \pi$ is called a Frobenius element at $p$. The kernel of the homomorphism $\pi_p\rTo \pi_p^u$ is denoted by $I_p$ and called the {\em inertia subgroup} at $p$.

The key interaction for applications to arithmetic are between $H^1(\Q, U)$ and the various $H^1(\Q_v, U)$. Since
$\pi_v$ injects into $\pi$, there is a restriction map
$$H^1(\Q,U)\rTo H^1(\Q_v,U)$$ for each $v$, which we put together into
$$\loc: H^1(\Q,U)\rTo \prod_v H^1(\Q_v, U).$$
Here then is the main problem of arithmetic gauge theory:
\begin{quote}
{\em For a gauge group $U$ over $\Q$, describe the image of $\loc$.}
\end{quote}
Any kind of a solution to this problem is called a {\em local-to-global principle} in number theory.

At this point, we pursue a bit more the analogy with geometric gauge fields. As discussed already, geometric gauge theory with symmetry group $U$ (in this case a real Lie group) deals with a space  $\cA$ of principal $U$-connections  on a spacetime manifold $X$. The usual convention these days is to take $\cA$ to be a space of $C^{\infty}$ connections. There is an action functional $$S:\cA\rTo \R$$
that is invariant under gauge transformations $\cU$ (connection preserving automorphisms of the principal bundle). The space of classical solutions is
$$M(X,U)=\cA^{EL}/\cU,$$
where $\cA^{EL}\subset \cA$ is the set of connections that satisfy the Euler-Lagrange equations for the functional $S$. The classical problem is to describe the space $M(X,U)$, or to find  points in $M(X,U)$ corresponding to specific boundary conditions. The quantum problem is to compute path integrals like
$$\int_{\cA/\cU} O_1(A)O_2(A)\cdots O_k(A) \exp(-S(A))dA,$$
where $O_i$ are local functions of $A$.

Now, from the point of view of classical physics, $M(X,U)$ will be the fields that we actually observe, and the embedding
$$M(X,U)\subset \cA/\cU$$
corresponds to a model for `quantum fluctuations' around classical solutions. However, the justification for considering $ \cA/\cU$ as the space of quantum fluctuation, or `off-shell states' of the field, is not so clear. It depends on the choice of an initial mathematical model inside which the classical solutions wre constructed. Some might argue that it is hard to even describe $M(X,U)$ unless one starts from $\cA/\cU$. This is false for physical reasons, since $M(X,U)$ is suppose to be a model of the classical states, which should make intrinsic sense regardless of the space in which we seek them\footnote{Of course this is not quite true. Classical states should be a statistical state of some sort arising out of the quantum theory, and hence, dependant on the quantum states. However, we are following here the tentative treatment found in standard expositions of path integral quantisation.}. Another objection comes from specific examples such as 3d Chern-Simons theory or 2d Yang-Mills theory, where $M(X,U)$ can easily be a space of {\em flat } connections. In that case, it has a topological description as a space of representations of the fundamental group of $X$. It is mostly this last case we have in mind when we consider the arithmetic versions. A model for quantum fluctuations of $M(X,U)$ might then just as well be collections of punctual local systems around  points of $X$, in the spirit of the jagged or singular paths that occur in Feynman's motivational description of the path integral \cite{zee}. This will be  especially appropriate  if we allow a model of $X$ that  can have complicated local topology.

It is from this point of view that we regard a collection $(P_v)_v$ of $\Q_v$ principal bundles as $v$ runs over the places of $\Q$ as a quantum arithmetic gauge field on $\Q$. A problem of finding   which collections glue together to a rational gauge field, that is, a principal $U$-bundle over $\Q$,  is the problem of describing the image of the localisation map. It is also an arithmetic analogue of finding and solving the Euler-Lagrange equation\footnote{This point of view was discovered independently by Philip Candelas and Xenia de la Ossa, although I may be misrepresenting their perspective.} . A better justification for this analogy will be discussed in section 10.

If we focus on a single component $P_v$, it is worth emphasising that the image of
$$\loc_v:  H^1(\Q,U)\rTo H^1(\Q_v, U)$$
consists exactly of those principal bundles $P_v$ whose external symmetry $\pi_v$ extends to the much larger group $\pi\supset \pi_v$. Computing this image precisely is critically related to the effective Mordell conjecture, as we will explain in section 8.
\section{The Tate-Shafarevich group and abelian gauge fields}

We should remark that the {\em kernel} of the localisation map is also frequently of importance. In words, these are the {\em locally trivial torsors\footnote{One of the complexities associated with arithmetic gauge fields is that many are not locally trivial in a naive sense, unlike the geometric situation. Part of the motivation for the \'etale topology is to have a topology that is fine enough so that natural torsors become locally trivial.}}. The best known case is when we have an elliptic curve $E$.
The kernel of localisation is then called the Tate-Shafarevich group of $E$, and denoted\footnote{Read  `Sha'.} $$\Sh(\Q, E).$$ A simple example  is when $E$ is given by the (non-Weierstrass) equation $$x^3+y^3+60z^3=0.$$ Then the curve $C$ given by $$3x^3+4y^3+5z^3=0$$ is an element\footnote{ This also is not so easy to see. There is an action of $E$ on $C$, which arises from the fact that $E$ is actually the Jacobian of $C$ \cite{AKMMMP}.} of $\Sh(\Q,E)$. The group $H^1(\Q,E)$ is an infinite torsion\footnote{To see this, one notes that all elements of $H^1(\Q,E)$ can be represented by the $\Qbar$ points of an algebraic curve $C$ in such a way that the action is algebraic and defined over $\Q$ \cite{silverman}. The torsor becomes trivial as soon as $C$ has a rational point. Now $C$ has a rational point over some finite field extension $K$ of $\Q$. That is, the class will become trivial under the restriction map $H^1(\Q, E)\rTo H^1(K,E)$. However, there is also a `trace' map $H^1(K,E)\rTo H^1(\Q,E)$ defined by summing a Galois conjugacy class of torsors. Also, the composed map $H^1(\Q, E)\rTo H^1(K,E)\rTo H^1(\Q, E)$ is simply multiplication by $[K:\Q]$. Thus, the element $[C]\in H^1(\Q,E)$ is killed by this degree.
} group. Remarkably, $\Sh(\Q,E)$ is conjectured to be finite. 

The relationship between the localisation map and $\Sh$ is a crucial part of the so-called `descent algorithm' for computing the points on an elliptic curve.
Recall that $E(\Q)$ is a finitely generated abelian group, that is,
$$E(\Q)\simeq \Z^r\times \mbox{finite abelian group}$$ and that its torsion subgroup is easy to compute \cite{silverman}. However, the rank is still a difficult quantity and the conjecture of Birch and Swinnerton-Dyer (BSD) is mainly concerned with the computation of $r$. The standard method at the moment is to look at $$E(\Q)/pE(\Q)$$ for some prime $p$, often $p=2$. If we know the structure of this group, it is elementary group theory to figure out the rank of $E(\Q)$, given that  we also know the torsion. The long exact sequence arising from
$$0\rTo E[p]\rTo E\rTo^p E\rTo 0$$
Gives
$$0\rTo E(\Q)/pE(\Q)\rInto H^1(\Q, E[p])\rTo^i H^1(\Q,E)[p]\rTo 0.$$
We have $$\Sh(\Q,E)[p]\subset H^1(\Q,E)[p].$$ Define the $p$-Selmer group $Sel(\Q, E[p])$ to be the inverse image of $\Sh(\Q,E)[p]$ under the map $i$, so that it fits into an exact sequence
$$0\rTo E(\Q)/pE(\Q)\rTo Sel(\Q, E[p])\rTo \Sh(\Q,E)[p]\rTo 0.$$
Described in words, $Sel(\Q, E[p])$ consists of the $E[p]$-torsors that become locally trivial when pushed out to $E$-torsors.

The key point is that the Selmer group is effectively computable, and this already gives us a bound on the Mordell-Weil group of $E$. This is then refined by way of  the diagram
$$\bd 0&\rTo & E(\Q)/p^nE(\Q)&\rTo &Sel(\Q, E[p^n])&\rTo &\Sh(\Q,E)[p^n]&\rTo &0 \\
&&\dTo & & \dTo & & \dTo & &\\
0&\rTo & E(\Q)/pE(\Q)&\rTo &Sel(\Q, E[p])&\rTo &\Sh(\Q,E)[p]&\rTo &0 \ed$$
for increasing values of $n$.
Provided $\Sh$ is finite, one can see that the image of $E(\Q)/pE(\Q)$ in $Sel(\Q,E[p])$ consists exactly of the elements that can be lifted to $Sel(Q,E[p^n])$ for all $n$. We get thereby, a cohomological expression for the group $E(\Q)/pE(\Q)$ that can be used to compute its structure precisely. The idea is
to compute the image $$Im(Sel(Q,E[p^n]))\subset Sel(Q,E[p])$$ for each $n$ and simultaneously compute the image of $E(\Q)_{\leq n}$ in $Sel(Q,E[p])$. Here, $E(\Q)_{\leq n}$ consists of the point in $E(\Q)$ of height\footnote{The height $h(x,y)$ of a point $(x,y)\in E(\Q)$ is defined as follows. Write
$(x,y)=(s/r, t/r)$ for coprime integers $s,t,r$. Then $h(x,y):=\log \sup \{|s|,|t|, |r|\}$. The height of the origin is defined to be zero. Clearly,  there are only finitely many rational $(x,y)$ of height $\leq n$ for any $n$. } $\leq n$. This is a finite set that can be effectively computed: just look at the finite set of pairs $(x,y)$ of height $\leq n$ and see which ones satisfy the equation for $E$. Thus, we have an inclusion
$$Im(E(\Q)_{\leq n})\subset Im(Sel(\Q,E[p^n])\subset Sel(\Q,E[p]).$$
Assuming $\Sh$ is finite, we get $$Im(E(\Q)_{\leq n})=Im(Sel(\Q,E[p^n])$$
for $n$ sufficiently large, at which point we can conclude that $$E(\Q)/pE(\Q)=Im(Sel(\Q,E[p^n]).$$
This is a conditional\footnote{In the sense that it terminates only if $\Sh$, or more precisely, its $p$-primary part $\Sh[p^{\infty}]$ is finite.} algorithm for computing the rank, which is used by all the existing computer packages. With slightly more care, it also gives a set of generators for the group $E(\Q)$. In this sense, we eventually arrive at a conditional algorithm for `completely determining' $E(\Q)$. An important part (some would argue the {\em most} important part) of BSD is to remove the `conditional' aspect.
\section{Non-abelian gauge fields and Diophantine geometry}

We keep to the conventions of the previous section and assume $U$ to be a $p$-adic Lie group over $\Q$. We now make the following further assumption: There is a finite set $S$ of places containing $p$ and $\infty$ such that for all $v\notin S$, the action of $\pi_v$ on $U$ factors through the quotient $\pi_v^u\simeq \Gal(\bar{\F}_v/\F_v).$  Another way of saying this is that the action of the inertia subgroup $I_v$ is trivial. We say that $U$ is {\em unramified} at $v$. Geometrically, this corresponds to having a family of groups on $\Spec(\Z_S)$, where $\Z_S\subset \Q$ is the ring of $S$-integers\footnote{ In terms of scheme theory, the set underlying $\Spec(\Z_S)$ is the open subset of $\Spec(\Z)$ obtained by removing the primes in $S$}, i.e., rational numbers whose denominators are divisible only by primes in $S$. We will assume that the torsors $P$ satisfy the same condition. In terms of  $\Spec(\Z)$, these can be thought of as connections having singularity\footnote{The geometry of schemes is organised in such a way that $\Z$  becomes a ring of functions on $\Spec(\Z)$.  It is easy to imagine that $\Z_S$ then becomes the ring of functions with restricted poles.}only along the primes in $S$. We will refer to these as $S$-integral $U$-torsors.  We will denote by $$H^1(\Z_S, U)$$ the isomorphism classes of $S$-integral $U$-torsors.
The reason for introducing this notion is that the $U$ and $P$ that arise in nature are $S$-integral for some $S$.  The condition of being unramified clearly makes sense even when $P$ is just a torsor over $\Q_v$. We  denote by $H^1_u(\Q_v, U)$ the isomorphism classes of unramified $U$ torsors over $\Q_v$.

One additional condition that $U$ and its torsors are required to satisfy is that of being {\em crystalline} at $p$, a technical condition about which we will be quite vague\footnote{This is an enormous subject in the study  of Galois representations and pedagogical references are easy to find with just the key word `crystalline representation.'  A comprehensive survey is in \cite{fontaine}.  In the non-abelian situation, a treatment is given in \cite{kim1}}. There is a big topological $\Q_p$-algebra $B_{cr}$ called the ring of $p$-adic periods and the torsors are required to trivialise over $B_{cr}$. This is a condition that comes from geometry and is closely related to $p$-adic Hodge theory. The point is that because $U$ is a $p$-adic Lie group, it will very rarely happen that the action is actually unramified at $p$. The crystalline condition captures smooth behaviour nevertheless. We denote by $H^1_f(\pi_p, U)$ the torsors over $\Q_p$ that are crystalline.

With these assumptions, we denote by
$$\prod'H^1(\Q_v, U)$$
the isomorphism classes of tuples $(P_v)_v$ where $P_v$ is a $U$-torsor over $\Q_v$ with the property that all but finitely many $P_v$ are unramified and such that $P_p$ is crystalline. For the global version, denote by $H^1_f(\Z_S, U)$ the $U$ torsors over $\Q$ that are unramified outside $S$  and crystalline at $p$.
Thus, we get a map
$$\loc: H^1_f(\Z_S, U)\rTo \prod'H^1(\Q_v, U),$$
whose image we would like to compute. 

 The main examples are 
\ms

(1) The constant group $U=GL_n(\Z_p)$ or other $p$-adic Lie groups with trivial $G$-action.

\ms
 In this case, from the earlier description in terms of cocycles, it is easy to see that a $U$-torsor is simply a representation
 $$\rho:G\rTo U.$$
 By our earlier assumption, this representation is required to be unramified outside $S$ and crystalline at $p$. We will return to this important case in the next section.
 \ms

(2) The $\Q_p$-pro-unipotent fundamental group \cite{deligne, kim2} $$U=\pi_1(\bar{V}, b)_{\Q_p}$$  of a smooth projective variety $V$ over $\Q$ equipped with a rational base-point $b\in V(\Q)$. We assume that $V$ extends to a smooth projective family over $\Z_{S\setminus p}$. An abstract definition  of $U$ can be given starting from the profinite \'etale fundamental group $\pi_1(\bV,b)$:  $\pi_1(\bar{V}, b)_{\Q_p}$ is the universal pro-unipotent group\footnote{An algebraic group is unipotent if it can be represented as a group of uppertriangular matrices with 1s on the diagonal. A pro-unipotent group is a projective limit of unipotent algebraic groups.} over $\Q_p$ admitting a continuous homomorphism $$\pi_1(\bar{V}, b)\rTo \pi_1(\bar{V}, b)_{\Q_p}$$
This is one of a number of `algebraic envelopes' of a group that have been important in both arithmetic and algebraic geometry \cite{amoros}. In spite of the difficulty of definition, it is substantially easier to work with than either the `bare' fundamental group or its profinite completion.

\ms

The important  and convenient fact is that $H^1_f(\Z_S, U)$ has the structure of a pro-algebraic scheme over $\Q_p$ \cite{kim1}. Among the constructions discussed so far, this is the closest to gauge-theoretic moduli spaces in physics and geometry. For another rational point $x$, one can also define  $$P(x)=\pi_1(\bar{V};b,x)_{\Q_p}:=[\pi_1(\bV;b,x)\times U]/\pi_1(\bV, b),$$ the $U$ torsor of pro-unipotent paths from $b$ to $x$. This construction gives us a map
$$V(\Q)\rTo H^1_f(\Z_S,U);$$
$$x\mapsto P(x)$$
that fits into a diagram
$$\bd V(\Q)&\rTo &V(\Q_p)\\
\dTo^A & & \dTo^{A_p }\\
H^1_f(\Z_S, U)& \rTo^{loc_p} & H^1_f(\Q_p, U)\ed$$
Even though the localisation map needs to be studied as a whole, because $U$ is a $p$-adic Lie group, it will usually happen that the component at $p$ is the most informative, and we will concentrate on this for now. (We will explain below the role of the adelic points.)
\begin{conj} \cite{BDKW} Suppose $V$ is a smooth projective curve of genus $\geq 2$. Then
$$A_p^{-1}(Im(\loc_p))=V(\Q).$$
\end{conj}

In essence, the conjecture is saying that the rational points can be recovered as the intersection between $p$-adic points and the space of $S$-integral torsors
inside the space of $p$-adic torsors.
A number of other diagrams are relevant to this discussion.
$$\bd H^1_f(\Z_S, U)& \rTo^{loc} & \prod' _vH^1(\Q_v, U)\\
&\rdTo &\dTo \\
& & H^1_f(\Q_p)\ed$$
The right vertical arrow is just the projection to the component at $p$. The image of the horizontal arrow should be computed by a reciprocity law \cite{kimrec1, kimrec2}, which we view as a preliminary version of the arithmetic Euler-Lagrange equations in that it specifies which collection of local torsors glue to a global torsor.

The diagram
$$\bd  X(\Q_p)& &  \\
 \dTo^{A_p} &\rdTo^{A^{DR}} & \\
 H^1(\Q_p, U)&\rTo^D & U^{DR}/F^0\ed$$
 is used to clarify he structure of the local moduli space $H^1_f(\Q_p, U)$ and to translate the Euler-Lagrange equations into   equations satisfied by the $p$-adic points.
 The last object $U^{DR}$ is the De Rham fundamental group \cite{deligne} endowed with a Hodge flitration $F^i$, which can be computed explicitly in such a way that  the map $A^{DR}$ is also described explicitly in terms of $p$-adic iterated integrals. From this point of view, computing the E-L equation is the main tool for finding the points $V(\Q)$ \cite{kim3, kim4, BD1, BD2}.
 
 To make this practical we use the lower central series $$U=U^1\supset U^2\supset U^3\supset \cdots,$$
 where $U^n=[U,U^{n-1}]$. We denote by $U_n=U/U^{n+1}$ the corresponding quotients, which are then finite-dimensional algebraic groups. All the diagrams above can be replaced by truncated versions, for example,
 $$\bd V(\Q)&\rTo &V(\Q_p)\\
\dTo^{A_n }& & \dTo^{A_{n,p} }\\
H^1_f(\Z_S, U_n)& \rTo^{loc_p} & H^1_f(\Q_p, U_n)\ed$$
These iteratively give equations for $V(\Q)$ depending on a reciprocity law for the image of $H^1_f(\Z_S, U_n)$ in $\prod'H^1(\Q_v, U_n)$.
 
 We illustrate this process with one example \cite{DW1, DW2}, which we take to be affine because it is easier to describe than the projective case. Let $V=\mathbb{P}^1\setminus \{0,1,\infty\}$. When we take $n=2$ and $S=\{\infty, 2,p\}$ the image of $H^1_f(\Z_S, U_2)$ in
 $$H^1_f(\Q_p, U_2)\simeq \A^3=\{(x,y,z)\}$$
 is described by the equation\footnote{The isomorphism between the local moduli space and the affine three-space is also a consequence of $p$-adic Hodge theory \cite{kim2}.}
 $$z-(1/2)xy=0.$$
 When translated back to points, this yields the consequence that the 2-integral points $V(\Z_{\{2\}})$ are in the zero set of the function
 $$D_2(z)=\ell_2(z)+(1/2)\log(z)\log(1-z).$$
 Here,  $\log(z)$ is the $p-$adic logarithm that is defined by the usual power series in a neighbourhood of 1 and then continued to all of $\Q \setminus \{ 0\}$ via additivity and the condition $\log(p)=0$. The $p$-adic $k$-logarithm is defined for $k\geq 2$ by
 $$\ell_k(z)=\sum_{n=1}^{\infty} z^n/n^k$$
 in a neighbourhood of zero and analytically continued to $V(\Q_p)$ using Coleman integration \cite{kim1}.
 
 When we use the defining equations for $H^1_f(\Z_S, U_4)$, we find that $V((\Z_{\{2\}})$ is killed by the additional equation
 $$
\z_p(3)\ell_4(z)+(8/7)[\log^32/24+\ell_4(1/2)/\log 2]\log(z) \ell_3(z)$$
$$+[(4/21)(\log^32/24+\ell_4(1/2)/\log 2)+\z_p(3)/24]\log^3(z)\log(1-z)=0.$$
Here, $\z_p(s)$ is the Kubota-Leopold $p$-adic zeta function.
\ms

It is worth noting that this method for finding rational points is a surprising confluence of three ingredients:

\ms

(1) The method of Chabauty \cite{coleman}, which was, in retrospect, the case of abelian gauge groups. One ends up using the $p$-adic logarithm on the Jacobian of the curve without considering cohomology at all;
\ms

(2) The descent method for finding points on elliptic curves \cite{silverman}. This again is another version of the abelian case, where one uses Galois cohomology and the Selmer group. As described earlier, this method is central to the conjecture of Birch and Swinnerton-Dyer.
\ms

(3) The geometry of arithmetic gauge fields.
 
 \ms
 
 Meanwhile, we should note that the ability to compute the full set of rational points will still rely on having a methodology that computes the image of the global moduli space quite precisely. We will return to this point in section 10. However, the hope that this should always be possible stems from the algebraicity of the localisation map
 $$\loc_p:H^1_f(\Z_S, U_n) \rTo^{loc_p}  H^1_f(\Q_p, U_n).$$
This implies that the image is a constructible set for the Zariski topology, which therefore admits a finite polynomial description.
 There is an increasing collection of examples for which the image has a precise enough description for the rational points to be computed completely \cite{BD1, BD2, DW1, DW2}. A spectacular recent result \cite{BDMTV} carries this out for the modular curve $$X_{s}(13)=X(13)/C^+_s(13),$$
 where $X(13)$ is the smooth projective model of the modular curve parametrising elliptic curves with full level 13 structure and $C^+_s(13)\subset GL_2(\F_{13})$ is the normaliser of a split Cartan subgroup. They prove the remarkable
 \begin{thm}[Balakrishnan, Dogra, M\"uller, Tuitman, Vonk] $X_s(13)$ has exactly 7 rational points consisting of 6 CM points and a cusp.
 \end{thm}
This theorem resolves a well-known difficulty in the arithmetic of modular curves arising in relation to effective versions of Serre's open image theorem. For a review of past work on this problem, see \cite{BP} and \cite{BPR}, where  13 is referred to as the `cursed' level.

\section{Galois representations, $L$-functions, and Chern-Simons actions}
We now consider $$H^1_f(\Z_S, GL_n(\Z_p)),$$ a moduli space of Galois representations\footnote{Here, we will allow ourselves to use  the word `space' quite loosely. There are numerous ways to geometrise this set, sometimes formally \cite{mazur3}, sometimes analytically \cite{chenevier}. It may also be most natural to regard it as a (derived) stack without worrying too much about representability. We will reprise this theme in the next section.}. This is the subspace of $$\Hom_{cont}(\pi_1(\Spec(\Z_S)), GL_n(\Z_p))/GL_n(\Z_p)$$
constrained by the crystalline condition\footnote{It is possible to be more general using more notions from $p$-adic Hodge theory} at $p$, where the $GL_n(\Z_p)$ is acting on the space of continuous homomorphisms by conjugation.
\ms

 It is believed then that the image can be characterised by an $L$-function via a reciprocity law that one is tempted to view as an arithmetic action principle of sorts. That is, for any collection
$(P_v)_v$ of local principal bundles with $P_p\in H^1_f(\Q_p, GL_n(\Z_p))$ crystalline and $P_v\in H^1_u(\Q_v,GL_n(\Z_p))$ for $v\notin S$, one looks at
the complex-valued product \cite{kiml}
$$L((P_v)_v, s):=\prod_{v\neq \infty} \frac{1}{\det(I-v^{-s}Fr_v|P_v^{I_v})},$$ 
which formally amalgamates the information associated to all `local' functions
$$(P_v)_v\mapsto \Tr((v^{-s}Fr_v)^n|P_v^{I_v})$$
as we run over places $v$ and natural numbers $n$.

Assuming a rather large number of standard conjectures in the theory of motives, some necessary conditions for $(P_v)_v$ to be in the image of the localisation map  of an irreducible representation $P$ are as follows.
\ms

(1) Each of the $\det(I-v^{-s}Fr_v|P^{I_v}_v )$ should be polynomials of $v^{-s}$ with integral coefficients for $v\notin S$; for any $v$, the coefficients should be algebraic. We use the algebraicity to regard the polynomial as complex-valued.

(2) There is an integer $w$ such that the absolute values of the eigenvalues of $Fr_v$ are $v^{w/2}$ for $v\notin S$.
This implies that the product converges absolutely for $Re(s)>w/2+1$.

(3)  $L((P_v)_v,s)$ has analytic continuation to all of $\C$ and satisfies a  functional equation of the form
$$L((P_v)_v,s)=ab^sL((P_v)_v,w+1-s),$$
for some rational numbers $a,b$. This function should have no poles unless $w$ is even, the $P_v$ are one-dimensional, and $Fr_v$ acts as $v^{w/2}$ for all but finitely many $v$.

\ms
 Roughly speaking, that these statements are necessary is summarised under the rubric of the Fontaine-Mazur conjecture \cite{FM} and the Hasse-Weil conjecture \cite{kiml}.
Even though it is not clear if a conjecture is stated in the literature, it appears to be commonly believed that these conditions should also characterise all $(P_v)_v$ that are in the image of the localisation map (cf. \cite{taylor}).

There is a sense in which $L((P_v)_v,s)$ should be related to an action. The complex number $s$ itself parametrises representations of a somewhat more general type, namely belonging to the idele class group of $F$. That is, for each place $v$ of $\Q$, there is a  normalised absolute value $\| \cdot \|_v$, which come together to form the norm character
$$\A_{\Q}^{\times}\rTo \C^{\times};$$
$$(a_v)_v\mapsto N((a_v)_v):=\prod_v \|a_v \|_v.$$
This character and its complex powers $N(\cdot)^{-s}$ factor through the idele class group $\A_{\Q}^{\times}/(\Q)^{\times}$ and the $L$ value is  a complex amplitude associated to $(P_v)_v$ twisted by $N(\cdot)^{-s}$.  The infinite product expansion will hold only for a region of $s$, so that the conjectured analytic continuation is supposed to involve a move from a kind of `decomposable range of the parameter' to one that is not. When the continuation is carried out, it turns out to be natural to view it as a section of a determinant line bundle, which is a function only in certain regions \cite{FK, kato}, creating an analogy with the wave functions of topological quantum field theory \cite{atiyah}. 

The question of finding natural action functionals on spaces of principal bundles appears to be important not just for unity of the theory, but because of the hope that it might lead to a more efficient approach to the gauge-theoretic Diophantine geometry  alluded to in the previous section. We will elaborate on this point in the next section.
While an action on torsors for $\pi_1(\bV,b)_{\Q_p}$ seems hard to define, there is an approach to   a Chern-Simons action  of Galois representations.  
To describe this, we lapse now into more geometric language and reproduce the discussion from \cite{kimcs, CKKPY}, which, in turn, is based on \cite{DW}. 

Let $X=\Spec(\cmcal{O}_F)$, the spectrum of the ring of integers in a number field $F$. We  assume that $F$ is totally imaginary. Denote by $\Gm$ the \'etale sheaf that associates to a scheme the units in the global sections of its coordinate ring. The topological fact underlying the functional is the canonical isomorphism (\cite[p. 538]{mazur}):
\[\label{eqn:*}
\inv: H^3(X, \Gm)\simeq \Q/\Z. \tag{$*$}
\]
This map is deduced from the `invariant' map of local class field theory \cite{NSW}. We will therefore use the same name for a  range of isomorphisms having the same essential nature, for example,
\[\label{eqn:**}
\inv:H^3(X, \Z_p(1))\simeq \Z_p, \tag{$**$}
\]
where $\Z_p(1)=\invlim_i \mu_{p^i}$, and  $\mu_n\subset \Gm$ is the sheaf of $n$-th roots of 1. The pro-sheaf $\Z_p(1)$ is a very familiar coefficient system for \'etale cohomology and (\ref{eqn:**}) is reminiscent of the fundamental class of a compact oriented three manifold for singular cohomology. Such an analogy was noted by Mazur around 50 years ago \cite{mazur2} and has been developed rather systematically by a number of mathematicians, notably, Masanori Morishita \cite{morishita}. Within this circle of ideas is included the analogy between knots and primes, whereby the map
\[
\Spec(\cmcal{O}_F/\mathfrak{P}_v)\rightarrowtail X
\]
from the residue field of a prime $\mathfrak{P}_v$ should be similar to the inclusion of a knot.  Let $F_v$ be the completion of $F$ at the prime $v$ and $\cmcal{O}_{F_v}$ its valuation ring. If one takes this analogy seriously, the map 
$$
\Spec(\cmcal{O}_{F_v})\to X,
$$ 
should be similar to the inclusion of a handle-body around the knot, whereas 
$$
\Spec(F_v)\to X
$$ 
resembles the inclusion of its boundary torus\footnote{It is not clear to us that the topology of the boundary should really be a torus. This is reasonable if one thinks of the ambient space as a three-manifold. On the other hand, perhaps it's possible to have a notion of a knot in a {\em homology three-manifold} that has an exotic tubular neighbourhood? In any case, M. Kapranov has pointed out that a better analogy is with a Klein bottle.}. Given a finite set $S$ of primes, we consider the scheme
\[
X^S:=\Spec(\cmcal{O}_F[1/S])=X\setminus \{\mathfrak{P}_v\}_{v\in S}.
\]
Since a link complement is homotopic to the complement of a tubular neighbourhood, the analogy is then forced on us between $X^S$ and a  three manifold with boundary given by a union of tori, one for each `knot' in $S$. These  are basic morphisms in $3$ dimensional topological quantum field theory \cite{atiyah}. From this perspective, the coefficient system $\Gm$ of the first isomorphism is analogous to the $S^1$-coefficient important in Chern-Simons theory \cite{witten, DW}. A more direct analogue of $\Gm$ is the sheaf $\cmcal{O}_M^{\times}$ of invertible analytic functions on a complex variety $M$. However, for compact K\"ahler manifolds, the comparison isomorphism 
$$
H^1(M, S^1)\simeq H^1(M, \cmcal{O}_M^{\times})_0,
$$
where the subscript refers to the line bundles with trivial topological Chern class, is a consequence of Hodge theory. This indicates that in the \'etale setting with no natural constant sheaf of $S^1$'s, the familiar $\Gm$ has a topological nature, and can be regarded as a substitute. \ms

We now move to the  definition of the arithmetic Chern-Simons action just for the simple case of a finite unramified Galois representation. Let
\[
\pi:=\pi_1(X, \fb),
\]  
be the profinite \'etale fundamental group of $X$,
where we take 
$$
\fb: \Spec(\ov F)\to X
$$ 
to be the geometric point coming from an algebraic closure of $F$. 
Assume now that the group $\mu_n({\ov F})$ of $n$-th roots of unity is in $F$ and fix a trivialisation 
$\zeta_n: \zmod n \simeq \mu_n$. 
This induces the isomorphism
$$
\inv: H^3(X, \zmod n) \simeq H^3(X, \mu_n) \simeq \frac{1}{n}\Z/\Z.
$$
Now let $A$ be a finite group with trivial $G_F$-action and fix a class $c\in H^3(A, \zmod n)$.
 For $$[\rho]\in H^1(\pi, A),$$ we get a class
$$
\rho^*(c)\in H^3(\pi, \Z/n\Z)
$$ 
that depends only on the isomorphism class $[\rho]$. Denoting by $\inv$ also the composed map 
\[
\xymatrix{
H^3(\pi,  \zmod n) \ar[r] & H^3(X, \zmod n) \ar[r]^-{\inv}_-{\simeq} & \frac{1}{n}\Z/\Z.
}
\]
We get thereby a function
\[
\xyv{0.15}
\xymatrix{
CS_c:H^1(\pi(X), A) \ar[r] & \nZ;\\
\hspace{12mm} [\rho] \hspace{2mm}\ar@{|->}[r] & \inv(\rho^*(c)).
}
\]
This is the basic and easy case of the classical Chern-Simons action in the arithmetic setting.
There is a natural generalisation to the case where ramification is allowed and where the representation has $p$-adic coefficients. It is related to natural invariants of algebraic number theory such as extensions of ideal class groups and $n$-th power residues symbols \cite{CKKPY, CKKPPY}. One might hope for such constructions to be related at once to $L$-functions and to Euler-Lagrange equations even for the unipotent fundamental groups of the previous section. Indeed, the approaches to the BSD conjecture that go via the `main conjecture of Iwasawa theory' take the view that Selmer groups should be annihilated by $L$-functions \cite{kato}. The reader might notice that the analytic equations defining integral points in the previous section actually indicate some connection to $L$-functions, but in a way that remains mysterious. 

We note also  that the Langlands reciprocity conjecture \cite{langlands} has as its goal the rewriting of arithmetic $L$-functions quite generally in terms of automorphic $L$-functions. In view of the striking work \cite{KW}, it seems reasonable to expect the geometry of arithmetic gauge fields to play a key role in importing  quantum field theoretic dualities to arithmetic geometry.

\section{Lagrangian Intersections}
In this section, we outline some detailed reasons to expect an arithmetic action principle, once again based on rather precise analogies with the theory of three-manifolds. Sufficiently rich geometric foundations\footnote{Currently, foundational work is in progress in collaboration with Kai Behrend and Yakov Kremnitzer. }
 underlying the constructions to follow should come from either rigid analytic geometry as in \cite{chenevier}  or derived versions \cite{BK}, \cite{PY}. For the purpose of this informal exposition, we will progress as though the necessary geometry is already in place, and simply use the natural properties we need. Of course, the reader should beware  the lack of rigorous foundations at the moment.
 \ms
 
Let $X$ be a compact oriented 3-manifold and let
$$X=X_1\cup_{\Sigma}X_2$$
be a Heegard splitting of $X$. Let $R$ be a sheaf of groups on $X$. (The precise nature of $R$ will be left vague for the purposes of this introduction.) Associated to this data, we have moduli spaces $$M(X_1,R), \ \ \ M(X_2, R),  \ \ M(\Sigma,R)$$ of principal $R$-bundles on $X_1$, $X_2$, and their common boundary $\Sigma$. There are also restriction maps
$$M(X_1,R)\rTo^{r_1}  M(\Sigma,R)\overset{r_2}{\lTo} M(X_2,R),$$
and geometric invariants are constructed out of the intersection of the images \cite{atiyah-weyl}.
Quite remarkably, this kind of intersection is of central interest in number theory as well.
\ms

We continue the discussion of the analogy mentioned in the previous section. If $K$ is an algebraic number  field and $\O_K$ its ring of integers, then $X=\Spec(\O_K)$ should be like a compact 3-manifold. If $v$ is a place of $K$ and $\cP_v$ the corresponding maximal ideal, the inclusion
$$\Spec(k_v)\rInto  Y$$
of the residue field $k_v=\O_K/\cP_v$ is supposed to be analogous to a knot. Let $K_v$ be the completion of $K$ at  $v$ and $\O_v$ its ring of integers.  One can then view the spectrum of the completion $Z_v=\Spec(\O_v)$ as a handle-body around the knot and $T_v=\Spec(K_v)$ as the complement of the knot inside the handle-body, which will therefore be homotopic to the boundary. The space $X^v=X\setminus \{\cP_v\}$, is the analogue to the knot complement, which is  homotopic to the complement of the interior of the handle-body. That is, we are viewing $X^v$ and  $Z_v$ as giving a Heegard splitting of the arithmetic 3-fold $X$ with common boundary $T_v$:
$$\bd & & T_v & & \\
&\ldTo & & \rdTo & \\
X^v & & & & Z_v \ed.$$ Now, given a sheaf
$R$ on $X$, we can again consider moduli spaces $M(X^v,R)$, $M(Z_v,R)$ and restriction maps
$$M(X^v,R)\rTo^{\loc_v} M(T_v,R)\overset{r_v}{\lTo} M(Z_v,R).$$
It turns out the intersection
 of these images are typically of great significance even in arithmetic. 
\ms

More generally, we can let $X^S=X\setminus \{\cP_{v}\}_{v\in S}$
for a finite set $S$ of places and consider maps
$$M(X^v,R)\rTo^{\loc_S}\prod_{v\in S} M(T_v,R)\overset{r_S}{\lTo} \prod_{v\in S}M(Z_v,R)$$
together with the corresponding fibre product
$$\mathcal{S}(X,R):=M(X^S, R)\times_{\prod_v M(T_v, R) }\prod_{v\in S} M(Z_v, R)$$
as well as other notions of intersection, naive or derived.
\ms

 We will be considering one of  two types.
\ms

\ms

(1) $R$ is a unipotent $\Q_p$-algebraic group with a continuous action of $\pi_S=\pi_1(X^S).$ In this case, we will assume that $S$ contains all places dividing $p$. Then the moduli space $M(X^S,R)$ can be identified with $$H^1(\pi_S, R)$$ representing the isomorphism classes of $\pi_S$-equivariant $R$-torsors. Similarly, $M(T_v, R)$ is the local cohomology $H^1(\pi_v, R)$, where $\pi_v=\Gal(\bF_v/F_v)$. The definition of $M(Z_v, R)$ is somewhat delicate. When $v$ doesn't divide $p$, then we let
$$M(Z_v, R):=H^1(\pi_v/I_v, R^{I_v}),$$
where $I_v\subset \pi_v$ is the inertia subgroup and $R^{I_v}$ refers to the part fixed by $I_v$. When $v|p$, then we let
$$M(Z_v, R):=H^1_f(\pi_v, Crys_v(R)),$$
where $Crys_v(R)$ is the maximal $\pi_v$-subgroup of $R$ which is crystalline, and the subscript $H^1_f$ refers to $\pi_v$-equivariant torsors for $Crys_v(R)$ which are themselves crystalline.
\ms

(2) $R=GL_n(E)$, where $E$ is a finite extension of $\Q_p$, the ring of integers in it, or a finite quotient of the ring of integers, considered as a constant sheaf on $X$. In this case, the moduli space $M(X^S,R)$ is $H^1(\pi_S, GL_n(E))$, the space of representations on free $E$-modules of rank $n$, and  $M(T_v, R)$ is of course just $H^1(\pi_v, GL_n(E))$, the space of representations of the local Galois group $\pi_v$. The definition of
$M(Z_v, R)$ again needs to distinguish between the case $v\nmid p$, when it consists of the unramified representations of $\pi_v$, and $v|p$, in which case it's made up of the isomorphism classes of crystalline representations.  This last notion is somewhat delicate in the case of integral representations, and we will be somewhat sloppy about the correct general notion.

\ms

We illustrate the centrality of these intersections in arithmetic geometry with some examples.

\ms

1. Let $A$ be an elliptic curve over $K$. Let $R$ be the sheaf associated to $V_p=T_p(A)\otimes \Q_p$ for some prime $p$ with the property that $A$ has good reduction at all places in $S_p=\{v \ | \ v|p\}$. Let $S$ be a finite set of places containing $S_p$ and the places of bad reduction for $A$. So we are regarding $V_p$ 
as a lisse sheaf of algebraic groups on $X^S$ that is simply pushed forward to $X$.
If  $v\in S\setminus S_p$, we have the unramified cohomology
$H^1_f(\pi_v, V_p):=H^1(\pi_v/I_v, V_p^{I_v})$ consisting of torsors that admit a reduction of structure group to the unramified subgroup of $V_p$. For $v\in S_p$,  $H^1_f(\pi_v, V_p)\subset H^1(\pi_v, V_p)$ is the subspace of crystalline torsors. The maps above become 
$$H^1(\pi_S, V_p)\rTo^{\loc_S} \prod_{v\in S}H^1(\pi_v, V_p)\lTo^{r_S}\prod_{v\in S}H^1_f(\pi_v, V_p).$$
In this case, we take the intersection to be the fiber product
$$Sel(A,\Q_p):=H^1(\pi_S, V_p)\times_{ \prod_{v\in S}H^1(\pi_v, V_p)} \prod_{v\in S}H^1_f(\pi_v, V_p).$$
\begin{conj}
$$\dim_{\Q_p} Sel(A,\Q_p)=\rank_{\Z} A(K).$$
\end{conj}
This is simply a reformulation of the conjecture on the finiteness of the $p$-part of the Tate-Shafarevich group of $A$.

\ms

2. Let $C/K$ be a hyperbolic curve, compact or affine. We assume given a regular $\O_K$-model, which we will not make explicit. This time, we take the group to be $U=\pi_1(\bar{C}, b)_{\Q_p}$, the $\Q_p$-pro-unipotent geometric \'etale fundamental group of $C$ as in section 8. The constructions can now be refined to
$$M(X^S, U)\rTo^{\loc_S} \prod_{v\in S}M(T_v, U)\lTo^{r_S}\prod_{v\in S}M(Z_v, U).$$

The applications to Diophantine geometry is now refined to the following diagram:

$$\bd
C(\O_K )&\rTo & Im(\loc_S)\cap Im(r_S) & \lTo & \prod_{v\in S} C(\O_v)  \\
  & & & & \\
\dTo& & \dTo& & \dTo \\ 
M(X^S, U)&\rTo^{\loc_S}&\prod_{v\in S}M(T_v, U) &\lTo^{r_S} &\prod_{v\in S}M(Z_v, U)
\ed$$
More precisely,  we have the implication
\bq
The projection $$Im(\loc_S)\cap Im(r_S) \rTo M(T_v, U)$$  is non-dense for some $v|p$.
\ms

 $\Rightarrow$ $C(\O_K)$ is finite.
\eq
The fact that the image of $C(\O_v)$ lies in $M(Z_v,U)$ is the main result of \cite{KT}.
\ms

3. The intersections for groups of type (2) above should be related to {\em arithmetic Casson invariants}. In fact, consider the case where $E/\Q_p$ is a finite extension and  denote by $H_v$ a filtration on
$(E\otimes_{\Q_p} K_v)^n$ for each $v\mid p$. Given a crystalline representation
$$\rho_v: \pi_v\rTo GL_n(E),$$ Fontaine's theory \cite{fontaine} associates to it a filtered $\phi-$module $$D(\rho_v)=(E^n\otimes B_{cr})^{\pi_v},$$ which is an  $E\otimes_{\Q_p} K_v$-module of rank $n$. If $v\mid p$, then $M(Z_v, R)^{H_v}$ denotes the crystalline representations $\rho_v$ of $\pi_v$ such that $D(\rho_v)$ has Hodge type $H_v$. If $v\nmid p$, then  $M(Z_v, R)^{H_v}:=M(Z_v, R)$. The Fontaine-Mazur conjecture proposes that
\bq
 {\em The locus of irreducible representations in the fiber product
$$M(X^S, R)\times_{\prod_{v\in S} M(T_v, R)}\prod_{v\in S} M(Z_v, R)^{H_v}$$
is finite.}
\eq
\ms

\ms

In considering Casson invariants in topology, it has been important to use the technique of realising subspaces as Lagrangian submanifolds of symplectic manifolds and then considering their intersections. This will happen for $\pi_1$-representations typically when the structure group is semi-simple, for example, $SL_n$. In the arithmetic setting, because of the the occurrence of Tate twists in duality, such a `self-dual' situation is harder to arrange. In the case (2), one could, for example, pass to representations of $$\pi_S^{\infty}\subset \pi_S,$$obtained by adjoining to the base field all $p$-power roots of 1, which will be classified by an infinite-dimensional space in general. We go on to outline now an alternative construction of Lagrangian intersections. For the remainder of this section, we will use also the notation of continuous group cohomology when it contributes to clarity.

\ms
Let $R$ be  a sheaf of $p$-adic analytic groups on $X^S=\Spec(\O_F[1/S])$ of types (1)  or (2) above and let $L$ be its Lie algebra. Assume $S$ contains all places dividing $p$ and that the action of $\pi_S=\pi_1(X^S)$ is crystalline at all places dividing $p$. Define $$T^*(1)R:=L^*(1)\rtimes R.$$ This is a twisted cotangent bundle of $R$. It's easy to see that the twisting still gives a well-defined group on $X^S$ by checking the compatibility with the action of $\pi_S$. Let
$$\tilde{c}\in H^1(\pi_S, T^*(1)R)$$ and $c\in H^1(\pi_S, R)$ be its image under the natural projection
$$H^1(\pi_S, T^*(1)R)\rTo H^1(\pi_S, R).$$ 
We note for later reference that this projection is split. Whenever we can geometrise these classifying spaces in a reasonable way, the tangent spaces will be computed as
$$T_{\tc}H^1(\pi_S, T^*(1)R)\simeq H^1(\pi_S, (L(c))^*(1)\times L(c))$$
$$\simeq H^1(\pi_S, (L(c))^*(1))\times H^1(\pi_S, L(c)).$$
where $L(c)$ is $L$ with the $\pi_S$-action twisted by the adjoint action of the cocycle $c$. This is because the adjoint action of $T^*(1)R$ on its tangent space factors through $R$. 
Similarly, if $v$ is a place of $F$, $\pi_v=\Gal(\bar{F_v}/F_v)$,  $\tc_v$ a cocycle of $\pi_v$ with values in $T^*(1)R$, and $c_v$ its projection to $R$, then
$$T_{\tc_v}H^1(\pi_v, T^*(1)R)\simeq H^1(\pi_v, (L(c_v))^*(1))\times H^1(\pi_v, L(c_v)).$$
Now,
$$H^1(\pi_v, (L(c_v))^*(1))\times H^1(\pi_v, L(c_v))\simeq
H^1(\pi_v, L(c_v))^*\times H^1(\pi_v, L(c_v))$$
by local Tate duality \cite{NSW}.
Hence, it carries a natural structure of a symplectic vector space, whose symplectic form $\omega_v$ is given by
$$\omega_v((\phi, c), (\phi',c'))=\langle \phi , c'\rangle- \langle \phi', c\rangle.$$
By summing over $v$, we get a symplectic structure on
$$\prod_{v\in S}[H^1(\pi_v, (L(c_v))^*(1))\times H^1(\pi_v, L(c_v))].$$
Even though the precise geometric foundation needs to be worked out,
$$\prod_{v\in S} M(T_v, T^*(1)R)=\prod_{v\in S}H^1(\pi_v, T^*(1)R)$$
should then have the structure of a analytic symplectic variety.
By Poitou-Tate duality \cite{NSW}, the image of
$$M(X^S, T^*(1)R)=H^1(\pi_S, T^*(1)R)$$
under the localisation
$$\loc_S: M(X^S, T^*(1)R) \rTo  \prod_{v\in S} M(T_v, T^*(1)R)$$
will then be a Lagrangian subvariety. If the cocycle $\tc_v$ is crystalline or unramified, then
$$H^1_f(\pi_v, (L(c_v))^*(1))\times H^1_f(\pi_v, L(c_v))$$
and
$$H^1(\pi_v/I_v, [(L(c_v))^*(1)]^{I_v})\times H^1(\pi_v/I_v,[ L(c_v)]^{I_v})$$
are  Lagrangian inside $$H^1(\pi_v, (L(c_v))^*(1))\times H^1(\pi_v, L(c_v)).$$  Hence,
$$\prod_{v\in S} M(Z_v, T^*(1)R)\subset \prod_{v\in S} M(T_v, T^*(1)R).$$
will also acquire the structure of  a Lagrangian subvariety. Using this, we can construct the Lagrangian intersection
$$\mathcal{S}(X, T^*(1)R):=M(X^S, T^*(1)R)\times_{\prod_{v\in S} M(T_v, T^*(1)R) }\prod_{v\in S} M(Z_v, T^*(1)R),$$
which possesses a split projection map to the intersection
$$\mathcal{S}(X, R)=M(X^S, R)\times_{\prod_{v\in S} M(T_v, R) }\prod_{v\in S} M(Z_v, R)$$
of interest.
\ms

A key point of this discussion is  the general theorem of \cite{BBJ}, which states that Lagrangian intersections are locally the critical loci of functions. Hence, in the Diophantine case where the relevant moduli spaces are schemes, at a Zariski local level, all the moduli spaces 
$$\mathcal{S}(X, T^*(1)R)$$ arising via the twisted cotangent construction can be described as the solution to Euler-Lagrange equations arising from a least action principle. From this point of view, the space of quantum off-shell fields is
$$\varinjlim_S \prod_{v\in S} M(T_v, T^*(1)R).$$
(Or possibly a restricted direct product.)
  The problem remains to give a natural construction of a global functional on this space from which equations of Euler-Lagrange type for both the Lagrangian intersection $\mathcal{S}(X, T^*(1)R)$ and the non-Lagrangian intersection $\mathcal{S}(X, R)$  might be extracted. In its absence, a possible interpretation is that the realisation of global moduli spaces as Lagrangian intersections itself should be viewed as an action principle, with an actual action being only locally defined. However, this view doesn't seem to lead to helpful computational tools, which might be considered the main goal in the Diophantine case.
\ms
\ms

    Here are some simple and concrete examples of this construction.

    \ms
        4. Consider the coefficient group $\Z_p^*$. Then the twisted cotangent bundle is
$$T^*(1)\Z_p^*=\Q_p(1)\times \Z_p^*.$$
Let $X=\Spec(\Z)$,  $S=\{p\}$, $X^S=\Spec(\Z)\setminus S$, and $\pi_S$, $\pi_p$, etc. as usual. In this case,
$$H^1(\pi_p, T^*(1)\Z_p^*)\simeq H^1(\pi_p, \Q_p(1))\times H^1(\pi_p, \Z^*_p) $$
$$\simeq (\hQp)\otimes \Q_p\times \Hom(\hQp, \Z^*_p) \simeq (\Z_p^*\times \hZ)\otimes \Q_p\times \Hom(\Z_p^*\times \hZ, \Z_p^*)$$
$$\simeq  \Q_p\times \Q_p\times W\times U,$$
where $W=\Hom(\Z_p^*, \Z_p^*)\simeq \Z_p^*$ is sometimes called the ($\Z_p$-points of the) {\em weight space}, and the space $U=\Hom(\hZ, \Z_p^*)\simeq \Z_p^*$ can be thought of as the unramified characters of $\pi_p$. Here we have used abelian local class field theory \cite{CF}. If we denote by
$H^1_f(\pi_p, T^*(1)\Z_p^*)$ the crystalline subspace, it's a fact \cite{loeffler}\footnote{This is a mathoverflow post, but with a short and complete proof.}that
$$H^1_f(\pi_p, \Z_p^*)\simeq \chi_p^{\Z}\times U\subset W\times U,$$
where $\chi_p$ is the identity map, also thought of as the $p$-adic cyclotomic character of $\pi_p$. Meanwhile \cite{bloch-kato},
$$H^1_f(\pi_p, \Q_p(1))\simeq \Q_p\times 0\subset \Q_p\times \Q_p.$$
So
$$H^1_f(\pi_p, T^*(1)\Z_p^*)\simeq \Q_p\times 0\times \chi_p^{\Z}\times U \subset  \Q_p\times \Q_p\times W\times U.$$
Now we examine the image of $$H^1(\pi_S, T^*(1)\Z_p^*)\rTo^{\loc_p }H^1(\pi_p, T^*(1)\Z_p^*).$$
The space $H^1(\pi_S, \Z_p^*)$ consists of characters that are unramified outside $p$, which then must factor through a $p$-adic power of the cyclotomic character $\chi_p$. (We will use the same notation for the global cyclotomic character and its restriction to $\pi_p$.) Thus, it can be identified via localisation with $W$ above. Meanwhile, $H^1(\pi_S, \Q_p(1))$ is generated by the image of the $p$-units, and hence, is generated (modulo torsion) by the image of $p$. This is just the subspace $0\times \Q_p$ of $H^1(\pi_p, \Q_p(1))$. Therefore, the image of $H^1(\pi_S, T^*(1)\Z_p^*)$ is just
$$0\times \Q_p\times W\times 0.$$
Hence, the Lagrangian intersection is
$$H^1(\pi_S, T^*(1)\Z_p^*)\cap H^1_f(\pi_p, T^*(1)\Z_p^*)\simeq 0\times 0\times \chi^{\Z}\times 0 .$$
Even though it's rather trivial, one interesting aspect  of this computation is that we get the same zero-dimensional intersection for
$$\mathcal{S}(\Spec(\Z), T^*(1)\Z_p^*)=H^1(\pi_S, T^*(1)\Z_p^*)\cap H^1_f(\pi_p, T^*(1)\Z_p^*)$$ and $$\mathcal{S}(\Spec(\Z), \Z_p^*)=H^1(\pi_S, \Z_p^*)\cap H^1_f(\pi_p, \Z_p^*).$$

\ms

\ms

5. Now let $X=\Spec(\O_F)$ for a totally complex number field $F$ and $S$ be the set of places in $F$ that lie over $p$.  Let $H$ be the class group of $F$. We have, again by local class field theory,
$$H^1(\pi_v, \Z_p^*)\simeq \Hom(\hF,\Z_p^*)\simeq \Hom(O_{F_v}^*, \Z_p^*)\times \Hom(\hZ, \Z_p^*)$$
$$= W_{F_v}\times U_{F_v},$$
where we write
$W_{F_v}= \Hom(O_{F_v}^*, \Z_p^*)$ and $U_{F_v}=\Hom(\hZ, \Z_p^*)$. This product decomposition of course relies on a choice of uniformiser, and implicitly, a choice of a totally ramified extension corresponding to it via Lubin-Tate theory.
On the other hand
$$H^1(\pi_v, \Q_p(1))\simeq (O_{F_v}^*\times \hZ)\otimes \Q_p\simeq \Q^{d_v}_p\times \Q_p,$$
where $d_v=[F_v:\Q_p].$ The subspace $H^1_f(\pi_v, \Q_p(1))$ is again identified with the first factor $\Q^{d_v}_p\times 0$. The crystalline subspace
$H^1_f(\pi_v, \Z_p^*)$ is identified \cite{loeffler} with
$$\chi_v^{\Z}\times U_{F_v},$$ where $\chi_v=\chi_p\circ N_v$ and $$N_v:\O_{F_v}^*\rTo \Z_p^*$$ is the norm. Thus,
$H^1_f(\pi_v, \Q_p(1)\times \Z_p^* )$ is identified with
$$\Q^{d_v}_p\times 0\times \chi_v^{\Z}\times U_{F_v}.$$
Consider now the image of $$H^1(\pi_S, \Q_p(1))\subset \prod_{v\in S} H(\pi_v,  \Q_p(1))=\Q_p^d\times \Q_p^s, $$ where $d=[F:\Q]$ and $s=|S|$.
This will be the $\Q_p$-span of the $S$-units $(\O_F[1/S])^*$  via the Kummer map. But we have an exact sequence
$$1\rTo \O_F^*\rTo (\O_F[1/S])^*\rTo \prod_{v\in S}\Z$$
given by the valuations at $v\in S$, and the last map has finite cokernel, since each maximal ideal $\cP_v$ corresponding to a place $v\in S$ has finite order in the ideal class group. So we find that the image of $H^1(\pi_S, \Q_p(1))$ is $E_p\times \Q_p^s$, where $E_p$ is the subspace of $\Q_p^d$ generated by the global units. (According to the Leopoldt conjecture, $E_p$ should have dimension $r_1+r_2-1$.)
\ms

As for $H^1(\pi_S, \Z_p^*)$ it's the same as $\Hom(A_S, \Z_p^*)$, where $A_S$ is the Galois group of the maximal abelian extension of $F$ unramified ourside $S$ \cite{CF}. 
This group fits into an exact sequence
$$0\rTo (\prod_{v\in S}\O_{F_v}^*)/\overline{Im(\O_F^*)}\rTo A_S\rTo H\rTo 0.$$
For the sake of simplicity, will will assume that this sequence is split, so that
$$\Hom(A_S, \Z_p^*)\simeq \Hom ((\prod_{v\in S}\O_{F_v}^*)/\overline{Im(\O_F^*)}, \Z_p^*)\times \Hom(H, \Z_p^*).$$
Thus, the image in
$\prod_{v\in S}[ W_{F_v}\times U_{F_v}]$ is
$$[\prod_{v\in S} W_{F_v}]^{glob} \times  \loc_S(\Hom(H, \Z_p^*) ),$$
where the superscript  denotes the set of products of local characters that vanish on the global units and $ \loc_p(\Hom(H, \Z_p^*) )$ refers to the restriction of global unramified characters to the decomposition groups in $S$. Therefore, the intersection for
$\mathcal{S}(X, T^*(1)\Z_p^*)$ is
$$E_p\times 0\times [\prod_{v\in S}\chi_v^{\Z}]^{glob}\times \loc_S(\Hom(H, \Z_p^*) ).$$

As before, the last two factors are discrete, but the first factor introduces a non-transverse contribution, which perhaps should be viewed as topologically trivial, in some sense. To spell it out again, the two subspaces of which it's the intersection are
$$\prod_{v\in S}H^1_f(\pi_v, \Q_p(1)\times \Z_p^* )=\Q_p^d\times 0\times \prod_{v\in S} \chi_v^{\Z}\times \prod_{v\in S} U_{F_v} $$
and
$$H^1(\pi_S, \Q_p(1)\times \Z_p^* )=E_p\times \Q_p^s\times [\prod_{v\in S} W_{F_v}]^{glob} \times  \loc_S(\Hom(H, \Z_p^*) ) $$
in $$\prod_{v\in S} H^1(\pi_v, \Q_p(1)\times \Z_p^* )=\prod_{v\in S} [ \Q_p^{d_v}\times \Q_p\times W_{F_v}\times U_{F_v}].$$

\ms

\ms

\ms

\ms

In any case, the examples above indicate that the two intersections $\mathcal{S}(X, T^*(1)R)$ and $\mathcal{S}(X, R)$
should have essentially similar structures. Also interesting is that the class group appears naturally in the computation of an `arithmetic Casson invariant'.

\ms

\ms

\pagebreak
\begin{flushleft}
{\bf Acknowledgements}
\end{flushleft}

I am extremely grateful to Yang-Hui He for the invitation to contribute this article. I am greatly indebted to Michael Atiyah, Philip Candeas, Xenia de la Ossa, Kai Behrend, Ted Chinburg, Tudor Dimofte, Richard Eager, Jeff Harvey, Kobi Kremnitzer, Tony Pantev, Ingmar Saberi,   Bertrand Toen,  and Johannes Walcher for many illuminating conversations.
Finally, it is a pleasure to acknowledge the kind interest and encouragement of Kevin Hartnett, whose friendly but persistent queries motivated me to write the last section.

\end{document}